\documentclass[oneside,english]{amsart}
\usepackage[T1]{fontenc}
\usepackage[latin9]{inputenc}
\usepackage{amssymb}
\usepackage{esint}

\makeatletter
\numberwithin{equation}{section} 
\numberwithin{figure}{section} 
  \theoremstyle{plain}
  \newtheorem{thm}{Theorem}[section]
  \theoremstyle{plain}
  \newtheorem{cor}[thm]{Corollary}
  \theoremstyle{plain}
  \newtheorem{lem}[thm]{Lemma}
  \theoremstyle{remark}
  \newtheorem*{acknowledgement*}{Acknowledgement}

\usepackage{babel}
\makeatother

\begin{document}

\title{Explicit Inversions of Certain Matrices I}

\author{Ruiming Zhang}

\begin{abstract}
In this note, we demonstrate a method to invert some Hankel matrices
explicitly by using the kernel polynomials for the related classical
orthogonal polynomials.
\end{abstract}


\curraddr{School of Mathematical Sciences\\
Guangxi Normal University\\
Guilin City, Guangxi 541004\\
P. R. China.}


\email{ruimingzhang@yahoo.com}

\maketitle

\section{Introduction}

In the theory of orthogonal polynomials, We could calculate the determinants
of some Hankel matrices once we know the three term recurrence relation
for the associated orthogonal polynomials and vice versa. It is well-known
that the kernel polynomials of the orthogonal polynomials encodes
important information about the associated Hankel matrices. These
matrices are generalizations of the Hilbert matrices. In this note
we present a method to invert some Hankel matrices associated with
classical polynomials by using the kernel polynomials. 

The following theorem very is a well known fact from the theory of
orthogonal polynomials:

\begin{thm}
\label{thm:1.1}Given a probability measure $\mu$ on $\mathbb{R}$
with a support of infinite many points. Let us consider the Hilbert
space of $\mu$-measurable functions \begin{equation}
\mathcal{X}:=\left\{ f(x)|\int|f(x)|^{2}d\mu(x)<\infty\right\} \label{eq:1.1}\end{equation}
with the inner product defined as\begin{equation}
(f,g):=\int f(x)\overline{g(x)}d\mu(x),\quad f,g\in\mathcal{X}.\label{eq:1.2}\end{equation}
Assume that $\left\{ w_{n}(x)\right\} _{n=0}^{\infty}$ is a sequence
of linearly independent functions in $\mathcal{X}$ with $w_{0}(x)=1$.
Let

\begin{equation}
\alpha_{jk}:=\int w_{j}(x)\overline{w_{k}(x)}d\mu(x),\quad j,k=0,1...,\label{eq:1.3}\end{equation}
 \begin{equation}
\Pi_{n}:=\left(\begin{array}{cccc}
\alpha_{00} & \alpha_{01} & \dots & \alpha_{0n}\\
\alpha_{10} & \alpha_{11} & \dots & \alpha_{1n}\\
\vdots & \vdots & \vdots & \vdots\\
\alpha_{n0} & \alpha_{n1} & \dots & \alpha_{nn}\end{array}\right),\quad n\in\mathbb{N}\cup\left\{ 0\right\} ,\label{eq:1.4}\end{equation}
 and

\begin{equation}
\Delta_{n}:=\det\Pi_{n},\quad n\in\mathbb{N}\cup\left\{ 0\right\} .\label{eq:1.5}\end{equation}
Then the $n$-th orthonormal function with positive coefficient in
$w_{n}(x)$ is given by the formula

\begin{equation}
p_{n}(x)=\frac{1}{\sqrt{\Delta_{n}\Delta_{n-1}}}\det\left(\begin{array}{ccccc}
\alpha_{00} & \alpha_{01} & \alpha_{02} & \dots & \alpha_{0n}\\
\alpha_{10} & \alpha_{11} & \alpha_{12} & \dots & \alpha_{1n}\\
\vdots & \vdots & \vdots & \ddots & \vdots\\
\alpha_{n0} & \alpha_{n1} & \alpha_{n2} & \dots & \alpha_{nn}\\
w_{0}(x) & w_{1}(x) & w_{2}(x) & \dots & w_{n}(x)\end{array}\right)\label{eq:1.6}\end{equation}
 for $n\in\mathbb{N}$ , with\begin{equation}
p_{0}(x)=w_{0}(x)=1.\label{eq:1.7}\end{equation}
 Furthermore, the coefficient of $p_{n}(x)$ in $w_{n}(x)$ is \begin{equation}
\gamma_{n}:=\sqrt{\frac{\Delta_{n-1}}{\Delta_{n}}}.\label{eq:1.8}\end{equation}

\end{thm}
\begin{proof}
The proof is the same as for the case $w_{n}(x)=x^{n}$, which could
be found in any orthogonal polynomials textbooks such as \cite{Szego}.
\end{proof}
\begin{cor}
\label{cor:1.1} For $n\in\mathbb{N}$, we have \begin{equation}
\Delta_{n}=\prod_{k=1}^{n}\frac{1}{\gamma_{n}^{2}}.\label{eq:1.9}\end{equation}

\end{cor}
\begin{proof}
This is a trivial consequence of \eqref{eq:1.7} and \eqref{eq:1.8}. 
\end{proof}
\begin{lem}
\label{lem:1.1} Let \begin{equation}
k_{n}(x,y):=\sum_{k=0}^{n}p_{k}(x)\overline{p_{k}(y)},\quad n\in\mathbb{N}\cup\left\{ 0\right\} .\label{eq:1.10}\end{equation}
 Then, for any $\pi(x)$ in the linear span of $\left\{ w_{k}(x)\right\} _{0}^{n}$
, we have

\begin{equation}
\int\pi(x)\overline{k_{n}(x,y)}d\mu(x)=\pi(y).\label{eq:1.11}\end{equation}
 
\end{lem}
\begin{proof}
To see \eqref{eq:1.11}, just expand $\pi(x)$ in $p_{k}(x),k=0,1\dots,n$. 
\end{proof}
\begin{lem}
\label{lem:1.2} For each $n\in\mathbb{N}\cup\left\{ 0\right\} $,
the function $k_{n}(x,y)$ satisfying \eqref{eq:1.11} is unique.
\end{lem}
\begin{proof}
Suppose there are two such functions $h_{n}(x,y)$ and $k_{n}(x,y)$
, then, 

\begin{align}
0 & <||h_{n}(\cdot,y)-k_{n}(\cdot,y)||^{2}\label{eq:1.12}\\
= & (h_{n}(\cdot,y)-k_{n}(\cdot,y),h_{n}(\cdot,y)-k_{n}(\cdot,y))\nonumber \\
= & (h_{n}(\cdot,y)-k_{n}(\cdot,y),h_{n}(\cdot,y))-(h_{n}(\cdot,y)-k_{n}(\cdot,y),k_{n}(\cdot,y))\nonumber \\
= & 0,\nonumber \end{align}
 which is a contradiction. 
\end{proof}
\begin{lem}
\label{lem:1.3} Let \begin{align}
(\beta_{jk})_{0\le j,k\le n}: & =\Pi_{n}^{-1},\quad n\in\mathbb{N}\cup\left\{ 0\right\} .\label{eq:1.13}\end{align}
Then,\begin{align}
k_{n}(x,y) & =\sum_{j,k=0}^{n}\beta_{jk}\overline{w_{j}(y)}w_{k}(x).\label{eq:1.14}\end{align}

\end{lem}
\begin{proof}
Let \begin{align}
f(x) & =\sum_{k=0}^{n}u_{k}w_{k}(x),\label{eq:1.15}\end{align}
 then,

\begin{align}
 & (f(\cdot),\sum_{j,k=0}^{n}\beta_{jk}\overline{w_{j}(y)}w_{k}(\cdot))\label{eq:1.16}\\
= & \sum_{m=0}^{n}u_{m}(w_{m}(\cdot),k(\cdot,y))\nonumber \\
= & \sum_{m=0}^{n}u_{m}\sum_{j,k=0}^{n}\overline{\beta_{jk}}w_{j}(y)(w_{m},w_{k})\nonumber \\
= & \sum_{m=0}^{n}u_{m}\sum_{j=0}^{n}w_{j}(y)\sum_{k=0}^{n}\overline{\beta_{jk}\alpha_{km}}\nonumber \\
= & f(y).\nonumber \end{align}
 By Lemma \ref{lem:1.2}, we have\begin{align}
k_{n}(x,y) & =\sum_{j,k=0}^{n}\beta_{jk}\overline{w_{j}(y)}w_{k}(x).\label{eq:1.17}\end{align}

\end{proof}
\begin{cor}
\label{cor:1.3}The kernel in Lemma \ref{lem:1.1} is also given by

\begin{align}
k_{n}(x,y) & =-\frac{1}{\Delta_{n}}\det\left(\begin{array}{ccccc}
0 & 1 & \overline{w_{1}(y)} & \cdots & \overline{w_{n}(y)}\\
1 & \alpha_{00} & \alpha_{01} & \dots & \alpha_{0n}\\
w_{1}(x) & \alpha_{10} & \alpha_{11} & \dots & \alpha_{1n}\\
\vdots & \vdots & \vdots & \ddots & \vdots\\
w_{n}(x) & \alpha_{n0} & \alpha_{n1} & \dots & \alpha_{nn}\end{array}\right)\label{eq:1.18}\end{align}
 for $n\in\mathbb{N}\cup\left\{ 0\right\} $. 
\end{cor}
\begin{proof}
Since \begin{align}
k_{n}(x,y) & =\sum_{j,k=0}^{n}\beta_{jk}\overline{w_{j}(y)}w_{k}(x),\label{eq:1.19}\end{align}
 with \begin{align}
\left(\beta_{jk}\right)_{0\le i,j\le n} & =\Pi_{n}^{-1}.\label{eq:1.20}\end{align}
 Then,\begin{align}
\beta_{jk} & =\frac{\Pi_{n}(k,j)}{\det\Pi_{n}}=\frac{\Pi_{n}(k,j)}{\Delta_{n}},\label{eq:1.21}\end{align}
 where $\Pi_{n}(k,j)$ is the $(k,j)$-th co-factor. Therefore,\begin{align}
k_{n}(x,y) & =\frac{1}{\Delta_{n}}\sum_{j,k=0}^{n}\Pi_{n}(k,j)\overline{w_{j}(y)}w_{k}(x).\label{eq:1.22}\end{align}
 It is clear that\begin{align}
\sum_{j,k=0}^{n}\Pi_{n}(k,j)\overline{w_{j}(y)}w_{k}(x) & =-\left|\begin{array}{cc}
0 & \mathbf{(\overline{W(y)})}^{T}\\
\mathbf{W(x)} & \Pi_{n}\end{array}\right|,\label{eq:1.23}\end{align}
 by direct determination expansion, which is \begin{align}
k_{n}(x,y) & =-\frac{1}{\Delta_{n}}\left|\begin{array}{cc}
0 & \mathbf{(\overline{W(y)})}^{T}\\
\mathbf{W(x)} & \Pi_{n}\end{array}\right|,\label{eq:1.24}\end{align}
where\begin{align}
\mathbf{W(x)} & =\left(\begin{array}{c}
1\\
w_{1}(x)\\
\vdots\\
w_{n}(x)\end{array}\right),\label{eq:1.25}\end{align}
and\begin{align}
\mathbf{\mathbf{(\overline{W(y)})}^{T}} & =\left(\begin{array}{cccc}
1, & \overline{w_{1}(y)} & ,\cdots, & \overline{w_{n}(y)}\end{array}\right).\label{eq:1.26}\end{align}

\end{proof}
Lemma \ref{lem:1.3} enables us to compute the inverse the Gram matrix
in terms of the orthonormal functions $\left\{ p_{n}(x)\right\} _{n=0}^{\infty}$.

\begin{cor}
\label{cor:1.4} Assume that $\left\{ w_{n}(x)\right\} _{n=0}^{\infty}$,
$\left\{ p_{n}(x)\right\} _{n=0}^{\infty}$ and $\Pi_{n}=(\alpha_{jk})_{0\le j,k\le n}$
as in Theorem \ref{thm:1.1}. Suppose we have two families of linear
functionals $\left\{ u_{k}\right\} _{k=0}^{\infty}$ and $\left\{ v_{k}\right\} _{k=0}^{\infty}$
over the linear space generated by $\left\{ w_{n}(x)\right\} _{n=0}^{\infty}$
with\begin{align}
u_{j}(w_{k}) & =\delta_{jk},\label{eq:1.27}\end{align}
 and

\begin{align}
v_{j}(\overline{w_{k}}) & =\delta_{jk}\label{eq:1.28}\end{align}
 for $j,k=0,1,...$. Then, 

\begin{align}
\beta_{jk} & =\sum_{m=0}^{n}u_{k}(p_{m}(x))v_{j}(\overline{p_{m}(y)}),\label{eq:1.29}\end{align}
where \begin{equation}
(\beta_{jk})_{0\le j,k\le n}=\Pi_{n}^{-1}.\label{eq:1.30}\end{equation}

\end{cor}
\begin{proof}
From Lemma \ref{lem:1.3}, we have\begin{equation}
\sum_{j,k=0}^{n}\beta_{jk}\overline{w_{j}(y)}w_{k}(x)=\sum_{m=0}^{n}\overline{p_{m}(y)}p_{m}(x).\label{eq:1.31}\end{equation}
Then we apply the functional $u_{j}$ and $v_{k}$ both sides of the
above equation, the claim of the corollary follows. 
\end{proof}

\section{Main Results\label{sec:Main-Results}}

\subsection{Preliminaries\label{sec:Preliminaries}}

The Euler's $\Gamma(z)$ is defined as \cite{Andrews} \begin{align}
\frac{1}{\Gamma(z)}: & =z\prod_{j=1}^{\infty}\left(1+\frac{z}{j}\right)\left(1+\frac{1}{j}\right)^{-z},\quad z\in\mathbb{C}\label{eq:2.1}\end{align}
For $a,a_{1},...,a_{r}\in\mathbb{C}$, the shifted factorials are
defined as\begin{align}
(a)_{n}: & =\frac{\Gamma(a+n)}{\Gamma(a)},\quad(a_{1},...,a_{r})_{n}:=\prod_{j=1}^{r}(a_{j})_{n},\quad n\in\mathbb{Z},r\in\mathbb{N}.\label{eq:2.2}\end{align}
The generalized hypergeometric series ${}_{r}F_{s}$ with parameters
$\left\{ a_{1},...,a_{r}\right\} $ and $\left\{ b_{1},...,b_{s}\right\} $
is formally defined by \begin{align}
{}_{r}F_{s}\left(\begin{array}{c}
a_{1},a_{2},...,a_{r}\\
b_{1},b_{2},...,b_{s}\end{array};z\right) & :=\sum_{n=0}^{\infty}\frac{(a_{1},...,a_{r})_{n}}{(a_{1},...,a_{s})_{n}}\frac{z^{n}}{n!}.\label{eq:2.3}\end{align}
The Barnes $G$-function is defined as\begin{align}
G(z): & =(2\pi)^{z/2}e^{-[z(z+1)+\gamma z^{2}]/2}\prod_{n=1}^{\infty}(1+\frac{z}{n})^{n}e^{-z+z^{2}/(2n)},\label{eq:2.4}\end{align}
where\begin{align}
\gamma: & =\lim_{n\to\infty}\left(\sum_{k=1}^{n}\frac{1}{k}-\ln n\right).\label{eq:2.5}\end{align}
The Barnes $G$-function is an entire function with the property\begin{align}
G(z+1) & =\Gamma(z)G(z),\label{eq:2.6}\end{align}
\begin{align}
\prod_{k=0}^{n}\Gamma(z+k) & =\frac{G(z+n+1)}{G(z)},\label{eq:2.7}\end{align}
and\begin{align}
G(n) & =\begin{cases}
0 & n=0,-1,-2,...\\
\prod_{i=0}^{n-2}i! & n=1,2,...\end{cases}.\label{eq:2.8}\end{align}
In following, in all the cases except the last one we use functionals\begin{align}
u_{i}(p(x)) & =v_{i}(p(x))=\frac{1}{i!}\left[\frac{d^{i}p(x)}{dx^{i}}\right]_{x=0},\label{eq:2.9}\end{align}
and for the last case we use \begin{align}
u_{i}(p(x)) & =v_{i}(p(x))=\frac{1}{i!}\left[\frac{d^{i}p(x)}{dx^{i}}\right]_{x=1},\label{eq:2.10}\end{align}
where $p(x)$ is a polynomial in variable $x$.

\subsection{The Hermite Polynomials }

The Hermite polynomials $\left\{ H_{n}(x)\right\} _{n=0}^{\infty}$
are defined as \cite{Andrews}\begin{align}
H_{n}(x) & =(2x)^{n}{}_{2}F_{0}\left(\begin{array}{c}
-\frac{n}{2},-\frac{n}{2}+\frac{1}{2}\\
-\end{array};-\frac{1}{x^{2}}\right)\label{eq:2.11}\end{align}
for $n\ge0$ and \begin{align}
H_{-1}(x) & =0.\label{eq:2.12}\end{align}
They satisfy\begin{equation}
DH_{n}(x)=2nH_{n-1}(x),\quad n\in\mathbb{N}\cup\left\{ 0\right\} .\label{eq:2.13}\end{equation}
Hermite polynomials satisfies 

\begin{align}
\int_{\mathbb{R}}H_{n}(x)H_{m}(x)\exp(-x^{2})dx & =2^{n}n!\sqrt{\pi}\delta_{mn}\label{eq:2.14}\end{align}
 for $n,m=0,1,...$. 

Thus, the orthonormal polynomials \begin{align}
h_{n}(x): & =\frac{H_{n}(x)}{\sqrt{n!2^{n}\sqrt{\pi}}}\label{eq:2.15}\end{align}
 have leading coefficients\begin{align}
\gamma_{n} & =\sqrt{\frac{2^{n}}{n!\sqrt{\pi}}}.\label{eq:2.16}\end{align}
Clearly,\begin{align}
\int_{-\infty}^{\infty}y^{n}e^{-y^{2}}dy & =\frac{1+(-1)^{n}}{2}\Gamma\left(\frac{n+1}{2}\right),\label{eq:2.17}\end{align}
 and

\begin{align}
\alpha_{ij} & =\frac{1+(-1)^{i+j}}{2}\Gamma\left(\frac{i+j+1}{2}\right),\label{eq:2.18}\end{align}
 for $i,j=0,1,...,n$. Thus, \begin{align}
\det\left(\frac{1+(-1)^{i+j}}{2}\Gamma\left(\frac{i+j+1}{2}\right)\right)_{j,k=0}^{n} & =2^{-\frac{n(n+1)}{2}}\pi^{\frac{n+1}{2}}\prod_{k=0}^{n}k!\label{eq:2.19}\end{align}
 or \begin{align}
\det\left(\frac{1+(-1)^{i+j}}{2}\Gamma\left(\frac{i+j+1}{2}\right)\right)_{j,k=0}^{n} & =2^{-\frac{n(n+1)}{2}}\pi^{\frac{n+1}{2}}G(n+2)\label{eq:2.20}\end{align}
 for $n=0,1,...$. 

The $(i,j)$-th entry of $\Pi_{n}^{-1}=(\beta_{jk})_{j,k=0}^{n}$
is\begin{align}
\beta_{ij} & =\frac{1}{i!j!}\sum_{k=\max(i,j)}^{n}\frac{1}{k!2^{k}\sqrt{\pi}}\left[\frac{d^{i}H_{k}(x)}{dx^{i}}\right]_{x=0}\left[\frac{d^{j}H_{k}(y)}{dy^{j}}\right]_{y=0}\label{eq:2.21}\\
= & \sum_{k=\max(i,j)}^{n}\frac{1}{k!2^{k}\sqrt{\pi}}\left\{ 2^{i}\binom{k}{i}H_{k-i}(0)\right\} \left\{ 2^{j}\binom{k}{j}H_{k-j}(0)\right\} \nonumber \end{align}
 or \begin{align}
\beta_{ij} & =\frac{2^{i+j}}{\sqrt{\pi}}\sum_{k=\max(i,j)}^{n}\frac{1}{k!2^{k}}\left\{ \binom{k}{i}H_{k-i}(0)\right\} \left\{ \binom{k}{j}H_{k-j}(0)\right\} .\label{eq:2.22}\end{align}
 for $i,j=0,1,...,n$.

\begin{thm}
For $n\in\mathbb{N}\cup\left\{ 0\right\} $, the matrix\begin{align}
 & \left(\frac{1+(-1)^{i+j}}{2\sqrt{\pi}}\Gamma\left(\frac{i+j+1}{2}\right)\right)_{0\le i,j\le n}\label{eq:2.23}\end{align}
 has the determinant\begin{align}
\det\left(\frac{1+(-1)^{i+j}}{2\sqrt{\pi}}\Gamma\left(\frac{i+j+1}{2}\right)\right)_{i,j=0}^{n} & =2^{-\frac{n(n+1)}{2}}G(n+2),\label{eq:2.24}\end{align}
 and its inverse matrix is \begin{align}
 & \left(\sum_{k=\max(i,j)}^{n}\frac{2^{i+j}\left\{ \binom{k}{i}H_{k-i}(0)\right\} \left\{ \binom{k}{j}H_{k-j}(0)\right\} }{k!2^{k}}\right)_{0\le i,j\le n}.\label{eq:2.25}\end{align}

\end{thm}

\subsubsection{The Laguerre Polynomials }

The Laguerre polynomials $\left\{ L_{n}^{\alpha}(x)\right\} _{n=0}^{\infty}$
may be defined as \cite{Andrews}\begin{align}
L_{n}^{\alpha}(x) & =\frac{(\alpha+1)_{n}}{n!}{}_{1}F_{1}\left(\begin{array}{c}
-n\\
\alpha+1\end{array};x\right)\label{eq:2.26}\end{align}
 for $n\ge0$, and we assume that\begin{align}
L_{-1}^{\alpha}(x) & =0.\label{eq:2.27}\end{align}
 We also have

\begin{align}
\frac{dL_{n}^{\alpha}(x)}{dx} & =-L_{n-1}^{\alpha+1}(x)\label{eq:2.28}\end{align}
and

\begin{align}
\int_{0}^{\infty}L_{m}^{\alpha}(x)L_{n}^{\alpha}(x)x^{\alpha}e^{-x}dx & =\frac{\Gamma(\alpha+n+1)}{n!}\delta_{mn},\label{eq:2.29}\end{align}
 for $\alpha>-1$ and $n,m=0,1,...$. The orthonormal polynomials

\begin{align}
p_{n}(x) & =(-1)^{n}\sqrt{\frac{n!}{\Gamma(\alpha+n+1)}}L_{n}^{\alpha}(x),\label{eq:2.30}\end{align}
 have the leading coefficients \begin{align}
\gamma_{n} & =\frac{1}{\sqrt{n!\Gamma(\alpha+n+1)}}\label{eq:2.31}\end{align}
 for $n=0,1,...$. Clearly

\begin{align}
\int_{0}^{\infty}x^{n+\alpha}e^{-x}dx & =\Gamma(\alpha+n+1),\label{eq:2.32}\end{align}
 and

\begin{align}
\alpha_{ij} & =\Gamma(\alpha+i+j+1)\label{eq:2.33}\end{align}
 for $i,j=0,1,...,n$, Then, \begin{align}
\det\left(\Gamma(\alpha+i+j+1)\right)_{j,k=0}^{n} & =\prod_{k=0}^{n}\left\{ k!\Gamma(\alpha+k+1)\right\} ,\label{eq:2.34}\end{align}
 or\begin{align}
\det\left(\Gamma(\alpha+i+j+1)\right)_{j,k=0}^{n} & =\frac{G(n+2)G(\alpha+n+2)}{G(\alpha+1)}.\label{eq:2.35}\end{align}
 Let $\Pi_{n}^{-1}=(\beta_{jk})_{j,k=0}^{n}$ , then, \begin{align}
\beta_{ij} & =\frac{1}{i!j!}\sum_{k=\max(i,j)}^{n}\frac{k!}{\Gamma(\alpha+k+1)}\left[\frac{d^{i}L_{k}^{\alpha}(x)}{dx^{i}}\right]_{x=0}\left[\frac{d^{j}L_{k}^{\alpha}(y)}{dy^{j}}\right]_{y=0}\label{eq:2.36}\\
= & \frac{1}{i!j!}\sum_{k=\max(i,j)}^{n}\frac{k!}{\Gamma(\alpha+k+1)}\left[(-1)^{i}L_{k-i}^{\alpha+i}(x)\right]_{x=0}\left[(-1)^{j}L_{k-j}^{\alpha+j}(y)\right]_{y=0},\nonumber \end{align}
 or \begin{align}
\beta_{ij} & =\frac{(-1)^{i+j}}{i!j!}\sum_{k=\max(i,j)}^{n}\frac{k!L_{k-i}^{\alpha+i}(0)L_{k-j}^{\alpha+j}(0)}{\Gamma(\alpha+k+1)}\label{eq:2.37}\end{align}
 for $j,k=0,1,...,n$. 

\begin{thm}
For $n=0,1,...$, the matrix\begin{align}
 & \left((\alpha+1)_{i+j}\right)_{0\le i,j\le n}\label{eq:2.38}\end{align}
 has determinant\begin{align}
\det\left((\alpha+1)_{i+j}\right)_{i,j=0}^{n} & =\frac{G(n+2)G(\alpha+n+2)}{G(\alpha+1)\Gamma(\alpha+1)^{n+1}},\label{eq:2.39}\end{align}
 and inverse \begin{align}
 & \left(\frac{\sum_{k=\max(i,j)}^{n}\frac{(\alpha+1)_{k}}{k!}\binom{k}{i}\binom{k}{j}}{(-1)^{i+j}(\alpha+1)_{i}(\alpha+1)_{j}}\right)_{0\le i,j\le n}\label{eq:2.40}\end{align}

\end{thm}

\subsubsection{The Ultraspherical Polynomials }

The Ultraspherical polynomials(or Genenbauer polynomials) $\left\{ C_{n}^{\lambda}(x)\right\} _{n=0}^{\infty}$
are defined as, \cite{Andrews}\begin{align}
C_{n}^{\lambda}(x) & =\frac{(2\lambda)_{n}}{n!}{}_{2}F_{1}\left(\begin{array}{c}
-n,2\lambda+n\\
\lambda+\frac{1}{2}\end{array};\frac{1-x}{2}\right)\label{eq:2.41}\end{align}
 for $n\ge0$, and we assume that\begin{align}
C_{-1}^{\lambda}(x) & =0.\label{eq:2.42}\end{align}
 We also have\begin{align}
\frac{dC_{n}^{\lambda}(x)}{dx} & =2\lambda C_{n-1}^{\lambda+1}(x),\label{eq:2.43}\end{align}

\begin{align}
\int_{-1}^{1}C_{m}^{\lambda}(x)C_{n}^{\lambda}(x)(1-x^{2})^{\lambda-\frac{1}{2}}dx & =\frac{\pi\Gamma(2\lambda+n)}{2^{2\lambda-1}n!(\lambda+n)[\Gamma(\lambda)]^{2}}\delta_{mn},\label{eq:2.44}\end{align}
 for $\lambda>-\frac{1}{2}$ and $n,m=0,1,...$. The orthonormal polynomials 

\begin{align}
p_{n}(x) & =\sqrt{\frac{2^{2\lambda-1}n!(\lambda+n)[\Gamma(\lambda)]^{2}}{\pi\Gamma(2\lambda+n)}}C_{n}^{\lambda}(x)\label{eq:2.45}\end{align}
have leading coefficients \begin{align}
\gamma_{n} & =\sqrt{\frac{(\lambda+n)2^{2\lambda+2n-1}\Gamma(\lambda+n)^{2}}{\pi n!\Gamma(2\lambda+n)}.}\label{eq:2.46}\end{align}
It is clear that

\begin{align}
\int_{-1}^{1}x^{n}(1-x^{2})^{\lambda-\frac{1}{2}}dx & =\frac{1+(-1)^{n}}{2}B\left(\frac{n+1}{2},\lambda+\frac{1}{2}\right),\label{eq:2.47}\end{align}
and

\begin{align}
\alpha_{ij} & =\frac{1+(-1)^{i+j}}{2}B\left(\frac{i+j+1}{2},\lambda+\frac{1}{2}\right),\label{eq:2.48}\end{align}
 for $i,j=0,1,...,n$, where $B(p,q)$ is the beta integral \begin{align}
B(p,q) & =\int_{0}^{1}x^{p-1}(1-x)^{q-1}dx,\quad\Re(p),\Re(q)>0.\label{eq:2.49}\end{align}
Then,\begin{align}
\det\left(\frac{1+(-1)^{i+j}}{2}B\left(\frac{i+j+1}{2},\lambda+\frac{1}{2}\right)\right) & =\prod_{k=0}^{n}\frac{\pi k!\Gamma(2\lambda+k)}{(\lambda+k)2^{2\lambda+2k-1}\Gamma(\lambda+k)^{2}},\label{eq:2.50}\end{align}
 or\begin{align}
\det\left(\frac{1+(-1)^{i+j}}{2}B\left(\frac{i+j+1}{2},\lambda+\frac{1}{2}\right)\right) & =\frac{\pi^{n+1}G(n+2)}{2^{(n+1)(n+2\lambda-1)}(\lambda)_{n+1}}\frac{G(2\lambda+n+1)G(\lambda)^{2}}{G(2\lambda)G(\lambda+n+1)^{2}}.\label{eq:2.51}\end{align}
 The $(i,j)$-th entry of the inverse matrix $\Pi_{n}^{-1}=(\beta_{jk})_{j,k=0}^{n}$
is\begin{align}
\beta_{ij} & =\frac{1}{i!j!}\sum_{k=\max(i,j)}^{n}\frac{2^{2\lambda-1}k!(\lambda+k)[\Gamma(\lambda)]^{2}}{\pi\Gamma(2\lambda+k)}\left[\frac{d^{i}C_{k}^{\lambda}(x)}{dx^{i}}\right]_{x=0}\left[\frac{d^{j}C_{k}^{\lambda}(y)}{dy^{j}}\right]_{y=0}\label{eq:2.52}\\
= & \frac{2^{2\lambda-1}[\Gamma(\lambda)]^{2}}{i!j!\pi}\sum_{k=\max(i,j)}^{n}\frac{k!(\lambda+k)}{\Gamma(2\lambda+k)}\left[2^{i}(\lambda)_{i}C_{k-i}^{\lambda+i}(x)\right]_{x=0}\left[2^{j}(\lambda)_{j}C_{k-j}^{\lambda+j}(x)\right]_{y=0},\nonumber \end{align}
 or \begin{align}
\beta_{ij} & =\frac{2^{i+j}(\lambda)_{i}(\lambda)_{j}\Gamma(\lambda)}{i!j!\sqrt{\pi}\Gamma(\lambda+\frac{1}{2})}\sum_{k=\max(i,j)}^{n}\frac{k!(\lambda+k)C_{k-i}^{\lambda+i}(0)C_{k-j}^{\lambda+j}(0)}{(2\lambda)_{k}}\label{eq:2.53}\end{align}
 for $i,j=0,1,...,n$.

\begin{thm}
Let $(\alpha_{ij})_{0\le i,j\le n}$ be the matrix with entries\begin{align}
\alpha_{ij} & =\frac{1+(-1)^{i+j}}{2}B\left(\frac{i+j+1}{2},\lambda+\frac{1}{2}\right)\label{eq:2.54}\end{align}
 for $i,j=0,1,...,n$, then, \begin{align}
\det\left(\frac{1+(-1)^{i+j}}{2}B\left(\frac{i+j+1}{2},\lambda+\frac{1}{2}\right)\right) & =\prod_{k=0}^{n}\frac{\pi k!\Gamma(2\lambda+k)}{(\lambda+k)2^{2\lambda+2k-1}\Gamma(\lambda+k)^{2}},\label{eq:2.55}\end{align}
 or\begin{align}
\det\left(\frac{1+(-1)^{i+j}}{2}B\left(\frac{i+j+1}{2},\lambda+\frac{1}{2}\right)\right) & =\frac{\pi^{n+1}G(n+2)}{2^{(n+1)(n+2\lambda-1)}(\lambda)_{n+1}}\frac{G(2\lambda+n+1)G(\lambda)^{2}}{G(2\lambda)G(\lambda+n+1)^{2}}.\label{eq:2.56}\end{align}
 The inverse matrix $(\beta_{ij})_{0\le i,j\le n}$ has entries\begin{align}
\beta_{ij} & =\frac{2^{i+j}(\lambda)_{i}(\lambda)_{j}\Gamma(\lambda)}{i!j!\sqrt{\pi}\Gamma(\lambda+\frac{1}{2})}\sum_{k=\max(i,j)}^{n}\frac{k!(\lambda+k)C_{k-i}^{\lambda+i}(0)C_{k-j}^{\lambda+j}(0)}{(2\lambda)_{k}}\label{eq:2.57}\end{align}
 for $i,j=0,1,...,n$. 
\end{thm}

\subsubsection{The Jacobi Polynomials }

The Jacobi polynomials $\left\{ P_{n}^{(\alpha,\beta)}(x)\right\} _{n=0}^{\infty}$
may be defined as \cite{Andrews,Szego}\begin{align}
P_{n}^{(\alpha,\beta)}(x) & =\frac{(\alpha+1)_{n}}{n!}{}_{2}F_{1}\left(\begin{array}{c}
-n;n+\alpha+\beta+1\\
\alpha+1\end{array};\frac{1-x}{2}\right)\label{eq:2.58}\end{align}
 for $n\ge0$, and \begin{align}
P_{-1}^{(\alpha,\beta)}(x) & =0.\label{eq:2.59}\end{align}
 We also have\begin{align}
\frac{dP_{n}^{(\alpha,\beta)}(x)}{dx} & =\frac{n+\alpha+\beta+1}{2}P_{n-1}^{(\alpha+1,\beta+1)}(x),\label{eq:2.60}\end{align}
 and

\begin{align}
\int_{-1}^{1}P_{m}^{(\alpha,\beta)}(x)P_{n}^{(\alpha,\beta)}(x)w(x)dx & =h_{n}\delta_{mn}\label{eq:2.61}\end{align}
 for $\alpha,\beta>-1$ and $n,m=0,1,...$ with\begin{align}
w(x): & =(1-x)^{\alpha}(1+x)^{\beta},\label{eq:2.62}\end{align}
 and\begin{align}
h_{n}: & =\frac{2^{\alpha+\beta+1}\Gamma(\alpha+n+1)\Gamma(\beta+n+1)}{(2n+\alpha+\beta+1)\Gamma(\alpha+\beta+n+1)n!}.\label{eq:2.63}\end{align}
 The orthonormal polynomials 

\begin{align}
p_{n}(x) & =\sqrt{\frac{(2n+\alpha+\beta+1)\Gamma(\alpha+\beta+n+1)n!}{2^{\alpha+\beta+1}\Gamma(\alpha+n+1)\Gamma(\beta+n+1)}}P_{n}^{(\alpha,\beta)}(x)\label{eq:2.64}\end{align}
 with leading coefficients\begin{align}
\gamma_{n} & =\frac{2^{n+(\alpha+\beta)/2}\Gamma\left(\frac{\alpha+\beta+1}{2}+n\right)\Gamma\left(\frac{\alpha+\beta+2}{2}+n\right)\sqrt{\frac{\alpha+\beta+1}{2}+n}}{\sqrt{n!\pi\Gamma(\alpha+\beta+n+1)\Gamma(\alpha+n+1)\Gamma(\beta+n+1)}}\label{eq:2.65}\end{align}
 for $n=0,1,...$. The moments of the Jacobi measure are

\begin{align}
\mu_{n} & =\frac{2^{\alpha+\beta+1}\Gamma(\alpha+1)\Gamma(\beta+1)}{(-1)^{n}\Gamma(\alpha+\beta+1)}{}_{2}F_{1}\left(\begin{array}{c}
-n,\beta+1\\
\alpha+\beta+1\end{array};2\right)\label{eq:2.66}\end{align}
 for $n=0,1,...$. Then $(i,j)$-th entry of the Hankel matrix $\Pi_{n}=(\alpha_{jk})_{j,k=0}^{n}$
is

\begin{align}
\alpha_{ij} & =\frac{(-1)^{i+j}\Gamma(\alpha+1)\Gamma(\beta+1)}{2^{-\alpha-\beta-1}\Gamma(\alpha+\beta+1)}{}_{2}F_{1}\left(\begin{array}{c}
-i-j,\beta+1\\
\alpha+\beta+1\end{array};2\right)\label{eq:2.67}\end{align}
 for $i,j=0,1,...,n$. Then,\begin{align}
\det\Pi_{n} & =\frac{\left(\frac{\pi}{2^{n+\alpha+\beta}}\right)^{n+1}}{\left(\frac{\alpha+\beta+1}{2}\right)_{n+1}}\frac{G(n+2)G\left(\frac{\alpha+\beta+1}{2}\right)^{2}G\left(\frac{\alpha+\beta+2}{2}\right)^{2}}{G\left(\frac{\alpha+\beta+3}{2}+n\right)^{2}G\left(\frac{\alpha+\beta+4}{2}+n\right)^{2}}\label{eq:68}\\
\times & \frac{G(\alpha+\beta+n+2)G(\alpha+n+2)G(\beta+n+2)}{G(\alpha+\beta+1)G(\alpha+1)G(\beta+1)},\nonumber \end{align}
 and the $(i,j)$-th entry of the inverse matrix $\Pi_{n}^{-1}$ is\begin{align}
\beta_{ij} & =\frac{1}{i!j!}\sum_{k=\max(i,j)}^{n}\frac{(2k+\alpha+\beta+1)\Gamma(\alpha+\beta+k+1)k!}{2^{\alpha+\beta+1}\Gamma(\alpha+k+1)\Gamma(\beta+k+1)}\label{eq:69}\\
\times & \left[\frac{d^{i}P_{k}^{(\alpha,\beta)}(x)}{dx^{i}}\right]_{x=0}\left[\frac{d^{j}P_{k}^{(\alpha,\beta)}(y)}{dy^{j}}\right]_{y=0}\nonumber \\
= & \frac{1}{i!j!}\sum_{k=\max(i,j)}^{n}\frac{(2k+\alpha+\beta+1)\Gamma(\alpha+\beta+k+1)k!}{2^{\alpha+\beta+1}\Gamma(\alpha+k+1)\Gamma(\beta+k+1)}\nonumber \\
\times & \left[\frac{(k+\alpha+\beta+1)_{i}P_{k-i}^{(\alpha+i,\beta+i)}(x)}{2^{i}}\right]_{x=0}\left[\frac{(k+\alpha+\beta+1)_{j}P_{k-j}^{(\alpha+j,\beta+j)}(y)}{2^{j}}\right]_{y=0},\nonumber \end{align}
 or \begin{align}
\beta_{ij} & =\sum_{k=\max(i,j)}^{n}\frac{(2k+\alpha+\beta+1)\Gamma(\alpha+\beta+k+1)k!}{\Gamma(\alpha+k+1)\Gamma(\beta+k+1)}\label{eq:70}\\
\times & \left\{ \frac{(k+\alpha+\beta+1)_{i}P_{k-i}^{(\alpha+i,\beta+i)}(0)\}\{(k+\alpha+\beta+1)_{j}P_{k-j}^{(\alpha+j,\beta+j)}(0)}{i!j!2^{\alpha+\beta+i+j+1}}\right\} .\nonumber \end{align}
 for $i,j=0,1,...,n$.

\begin{thm}
For $n=0,1,...$, the matrix\begin{align}
 & \left({}_{2}F_{1}\left(\begin{array}{c}
-i-j,\beta+1\\
\alpha+\beta+1\end{array};2\right)\right)_{0\le i,j\le n}\label{eq:2.71}\end{align}
 has the determinant\begin{align}
 & \det\left({}_{2}F_{1}\left(\begin{array}{c}
-i-j,\beta+1\\
\alpha+\beta+1\end{array};2\right)\right)_{j,k=0}^{n}\label{eq:2.72}\\
= & \left(\frac{\Gamma(\alpha+\beta+1)2^{-(2\alpha+2\beta+n+1)}\pi}{\Gamma(\alpha+\beta+1)\Gamma(\alpha+\beta+1)}\right)^{n+1}\nonumber \\
\times & \frac{G(n+2)G\left(\frac{\alpha+\beta+1}{2}\right)^{2}G\left(\frac{\alpha+\beta+2}{2}\right)^{2}}{G\left(\frac{\alpha+\beta+3}{2}+n\right)^{2}G\left(\frac{\alpha+\beta+4}{2}+n\right)^{2}}\nonumber \\
\times & \frac{G(\alpha+\beta+n+2)G(\alpha+n+2)G(\beta+n+2)}{\left(\frac{\alpha+\beta+1}{2}\right)_{n+1}G(\alpha+\beta+1)G(\alpha+1)G(\beta+1)},\nonumber \end{align}
 and its inverse matrix $\left(\beta_{ij}\right)_{0\le i,j\le n}$with\begin{align}
\beta_{ij} & =\sum_{k=\max(i,j)}^{n}\frac{k!(2k+\alpha+\beta+1)(\alpha+\beta+1)_{k}}{(-2)^{i+j}i!j!(\alpha+1)_{k}(\beta+1)_{k}}\label{eq:2.73}\\
\times & \left\{ (k+\alpha+\beta+1)_{i}P_{k-i}^{(\alpha+i,\beta+i)}(0)\}\{(k+\alpha+\beta+1)_{j}P_{k-j}^{(\alpha+j,\beta+j)}(0)\right\} .\nonumber \end{align}
 for $i,j=0,1,...,n$. 
\end{thm}
If we take the polynomial sequence\begin{align}
w_{n}(x) & =(x-1)^{n}\label{eq:2.74}\end{align}
for $n=0,1,...$, and the linear functionals defined in \eqref{eq:2.10}.
Then the $(i,j)$-th entry of $\Pi_{n}=(\alpha_{jk})_{j,k=0}^{n}$
is

\begin{align}
\alpha_{ij} & =\int_{-1}^{1}(x-1)^{i+j}w(x)dx\label{eq:2.75}\end{align}
 or

\begin{align}
\alpha_{ij} & =\frac{2^{\alpha+\beta+i+j+1}\Gamma(\alpha+i+j+1)\Gamma(\beta+1)}{(-1)^{i+j}\Gamma(\alpha+\beta+i+j+2)}\label{eq:2.76}\end{align}
 for $i,j=0,1,...,n$, and its determinant is given by\begin{align}
\det\Pi_{n} & =\frac{G(\alpha+n+1)G(\beta+n+1)G(\alpha+\beta+n+1)}{2^{(n+\alpha+\beta)(n+1)}G(\alpha+1)G(\beta+1)G(\alpha+\beta+1)}\label{eq:2.77}\\
\times & \frac{\pi^{n+1}G(n+2)G\left(\frac{\alpha+\beta+1}{2}\right)^{2}G\left(\frac{\alpha+\beta+2}{2}\right)^{2}}{\left(\frac{\alpha+\beta+1}{2}\right)_{n+1}G\left(\frac{\alpha+\beta+1}{2}+n+1\right)^{2}G\left(\frac{\alpha+\beta+2}{2}+n+1\right)^{2}},\nonumber \end{align}
 and the it inverse matrix has entries\begin{align}
\beta_{ij} & =\frac{1}{i!j!}\sum_{k=\max(i,j)}^{n}\frac{(2k+\alpha+\beta+1)\Gamma(\alpha+\beta+k+1)k!}{2^{\alpha+\beta+1}\Gamma(\alpha+k+1)\Gamma(\beta+k+1)}\label{eq:2.78}\\
\times & \left[\frac{d^{i}P_{k}^{(\alpha,\beta)}(x)}{dx^{i}}\right]_{x=1}\left[\frac{d^{j}P_{k}^{(\alpha,\beta)}(y)}{dy^{j}}\right]_{y=1}\nonumber \\
= & \sum_{k=\max(i,j)}^{n}\frac{(2k+\alpha+\beta+1)\Gamma(\alpha+\beta+k+1)k!}{2^{\alpha+\beta+1}\Gamma(\alpha+k+1)\Gamma(\beta+k+1)}\nonumber \\
\times & \left[\frac{(k+\alpha+\beta+1)_{i}P_{k-i}^{(\alpha+i,\beta+i)}(x)}{i!2^{i}}\right]_{x=1}\left[\frac{(k+\alpha+\beta+1)_{j}P_{k-j}^{(\alpha+j,\beta+j)}(y)}{j!2^{j}}\right]_{y=1},\nonumber \end{align}
 or \begin{align}
\beta_{ij} & =\sum_{k=\max(i,j)}^{n}\frac{(2k+\alpha+\beta+1)\Gamma(\alpha+\beta+k+1)k!}{\Gamma(\alpha+k+1)\Gamma(\beta+k+1)}\label{eq:2.79}\\
\times & \left\{ \frac{(k+\alpha+\beta+1)_{i}P_{k-i}^{(\alpha+i,\beta+i)}(1)\}\{(k+\alpha+\beta+1)_{j}P_{k-j}^{(\alpha+j,\beta+j)}(1)}{i!j!2^{\alpha+\beta+i+j+1}}\right\} \nonumber \end{align}
 for $i,j=0,1,...,n$. Since

\begin{align}
P_{n}^{(\alpha,\beta)}(1) & =\frac{(\alpha+1)_{n}}{n!},\label{eq:2.80}\end{align}
 then

\begin{align}
\beta_{ij} & =\frac{\Gamma(\alpha+\beta+1)(\alpha+\beta+1)_{i}(\alpha+\beta+1)_{j}}{2^{\alpha+\beta+i+j+1}(\alpha+1)_{i}(\alpha+1)_{j}\Gamma(\alpha+1)\Gamma(\beta+1)}\label{eq:2.81}\\
\times & \sum_{k=\max(i,j)}^{n}\left\{ \frac{(2k+\alpha+\beta+1)(\alpha+1)_{k}}{k!(\alpha+\beta+1)_{k}(\beta+1)_{k}}\right\} \nonumber \\
\times & \left\{ \binom{k}{i}\binom{k}{j}(\alpha+\beta+i+1)_{k}(\alpha+\beta+j+1)_{k}\right\} \nonumber \end{align}
 for $i,j=0,1,...,n$. Therefore, we have proved the following:

\begin{thm}
For $n=0,1,...$, the determinant of the matrix \begin{align}
 & \left(\frac{(\alpha+1)_{i+j}}{(\alpha+\beta+2)_{i+j}}\right)_{0\le i,j\le n}\label{eq:2.82}\end{align}
 is \begin{align}
 & \det\left(\frac{(\alpha+1)_{i+j}}{(\alpha+\beta+2)_{i+j}}\right)_{i,j=0}^{n}\label{eq:2.83}\\
= & \frac{G\left(\frac{\alpha+\beta+1}{2}\right)^{2}G\left(\frac{\alpha+\beta+2}{2}\right)^{2}}{G(\alpha+1)G(\beta+1)G(\alpha+\beta+1)}\nonumber \\
\times & \left(\frac{\pi\Gamma(\alpha+\beta+2)}{2^{2n+2\alpha+2\beta+1}\Gamma(\alpha+1)\Gamma(\beta+1)}\right)^{n+1}\nonumber \\
\times & \frac{G(n+2)G(\alpha+n+1)G(\beta+n+1)G(\alpha+\beta+n+1)}{\left(\frac{\alpha+\beta+1}{2}\right)_{n+1}G\left(\frac{\alpha+\beta+1}{2}+n+1\right)^{2}G\left(\frac{\alpha+\beta+2}{2}+n+1\right)^{2}},\nonumber \end{align}
 and its inverse matrix $\left(\gamma_{ij}\right)_{0\le i,j\le n}$
has elements\begin{align}
\gamma_{ij} & =\frac{(-1)^{i+j}(\alpha+\beta+1)_{i}(\alpha+\beta+1)_{j}}{(\alpha+1)_{i}(\alpha+1)_{j}(\alpha+\beta+1)}\label{eq:2.83}\\
\times & \sum_{k=\max(i,j)}^{n}\left\{ \frac{(2k+\alpha+\beta+1)(\alpha+1)_{k}}{k!(\alpha+\beta+1)_{k}(\beta+1)_{k}}\right\} \nonumber \\
\times & \left\{ \binom{k}{i}\binom{k}{j}(\alpha+\beta+i+1)_{k}(\alpha+\beta+j+1)_{k}\right\} \nonumber \end{align}
 for $i,j=0,1,...,n$. 
\end{thm}
\begin{acknowledgement*}
This work is partially supported by Chinese National Natural Science
Foundation grant No.10761002, Guangxi Natural Science Foundation grant
No.0728090.
\end{acknowledgement*}

\end{document}